\documentclass[10pt]{amsart}
\usepackage{amsmath,amsthm,amssymb,amscd}
\usepackage{epsfig}
\usepackage{eucal}
\begin{document}
\newcommand{\con} {{\rm config}}
\newcommand{\Con} {{\rm Config}}
\newcommand{\Ba} [1] {\ensuremath{{\mathcal {B}}^{#1}}}    
\newcommand{\Hu} {{\mathcal H}}                            
\newcommand{\ba}{{\mathfrak b}}                            
\newcommand{\sba}{{\mathfrak {sb}}}                        
\newcommand{\pba}{{\mathfrak {pb}}}                        
\newcommand{\D} {{\mathfrak D}}                            
\newcommand{\B} {{\bf B}}                                  
\newcommand{\M}  [1] {\ensuremath{{\overline{\mathcal M}}{^{#1}_0(\R)}}}   
\newcommand{\oM} [1] {\ensuremath{{\mathcal M}_{0}^{#1}(\R)}}              
\newcommand{\dM} [1] {\ensuremath{{\widetilde{\mathcal M}}_{0}^{#1}(\R)}}  
\newcommand{\Sh} {{\mathcal{S}}}                           
\newcommand{\Pb} {{\bf P}}                                 
\newcommand{\C} {{\mathbb C}}                              
\newcommand{\R} {{\mathbb R}}                              
\newcommand{\Z} {{\mathbb Z}}                              
\newcommand{\Pj} {{\mathbb P}}                             
\newcommand{\T} {{\mathbb T}}                              
\newcommand{\Sg} {\mathbb S}                               
\newcommand{\sg} {\sigma}                                  
\newcommand{\Gl} {{\rm Gl}}                                
\newcommand{\G} {{\mathcal G}}                             
\newcommand{\Gc} {{\mathfrak S}}                           
\newcommand{\Cox} [1] {\ensuremath{J_{#1}}}                
\newcommand{\qCox} [1] {\ensuremath{\tilde{J}_{#1}}}       
\newcommand{\Bnd} [2] {{\ensuremath{\Lambda (#1,#2)}}}
\newcommand{\Cobnd} [2] {{\ensuremath{\Lambda^n(#1,#2)}}}
\newcommand{\Op} [1] {{\mathcal {O}}(#1)}                  
\newcommand{\Lc} [1] {{\mathcal {C}}(#1)}                  
\newcommand{\Li} [1] {{\mathcal {I}}(#1)}                  
\newcommand{\SI} {\ensuremath{SI}}                         

\theoremstyle{plain}
\newtheorem{thm}{Theorem}[subsection]
\newtheorem{prop}[thm]{Proposition}
\newtheorem{cor}[thm]{Corollary}
\newtheorem{lem}[thm]{Lemma}
\newtheorem{conj}[thm]{Conjecture}

\theoremstyle{definition}
\newtheorem{defn}[thm]{Definition}
\newtheorem{exmp}[thm]{Example}

\theoremstyle{remark}
\newtheorem*{rem}{Remark}
\newtheorem*{hnote}{Historical Note}
\newtheorem*{nota}{Notation}
\newtheorem*{ack}{Acknowledgments}

\numberwithin{equation}{section}

\title {Tessellations of moduli spaces and the mosaic operad}

\author{Satyan L. Devadoss}
\address{Department of Mathematics, Johns Hopkins University, Baltimore, Maryland 21218}
\email{devadoss@math.jhu.edu}

\begin{abstract}
We construct a new (cyclic) operad of \emph{mosaics} defined by polygons with marked diagonals. Its underlying (aspherical) spaces are the sets \M{n} which are naturally tiled by Stasheff associahedra.  We describe them as iterated blow-ups and show that their fundamental groups form an operad with similarities to the operad of braid groups.
\end {abstract}

\maketitle

{\small \begin{ack}
This paper is a version of my doctorate thesis under Jack Morava, to whom I am indebted for providing much guidance and encouragement.  Work of Davis, Januszkiewicz, and Scott has motivated this project from the beginning and I would like to thank them for many useful insights and discussions.  A letter from Professor Hirzebruch also provided inspiration at an early stage.  I am especially grateful to Jim Stasheff for bringing up numerous questions and for his continuing enthusiasm about this work.
\end{ack}}


\baselineskip=15pt

%
%
\section {The Operads}

\subsection{} \label{operad}
The notion of an operad was created for the study of iterated loop spaces~\cite{may1}.  Since then, operads have been used as universal objects representing a wide range of algebraic concepts.  We give a brief definition and provide classic examples to highlight the issues to be discussed.

\begin{defn}
An {\em operad} $\{\Op{n} \; | \; n \in \mathbb {N} \}$ is a collection of objects $\Op{n}$ in a monoidal category endowed with certain extra structures: 

1. $\Op{n}$ carries an action of the symmetric group $\Sg_n$.

2. There are composition maps
\begin{equation}
\Op{n} \otimes \Op{k_1} \otimes \cdots \otimes \Op{k_n} \rightarrow \Op{k_1 + \cdots + k_n}
\label{e:operad}
\end{equation}
\indent \hspace{10pt} which satisfy certain well-known axioms, {\em cf}.\ \cite{may2}.
\end{defn}

This paper will be concerned mostly with operads in the context of topological spaces, where the objects $\Op{n}$ will be equivalence classes of geometric objects.

\begin{figure} [h]
\centering {\includegraphics {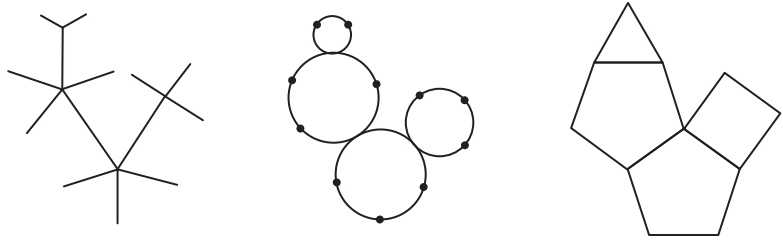}}
\caption{Trees, Bubbles, and Polygons}
\label{btp}
\end{figure}

\begin{exmp}
These objects can be pictured as {\em trees} (Figure~\ref{btp}a). A tree is composed of corollas\footnote{A corolla is a collection of edges meeting at a common vertex.} with one external edge marked as a {\em root} and the remaining external edges as {\em leaves}.  Given trees $s$ and $t$, basic compositions are  defined as $s \circ_i t$, obtained by grafting the root of $s$ to the $i^{\rm th}$ leaf of $t$. This grafted piece of the tree is called a {\em branch}.
\end{exmp}

\begin{exmp}
There is a dual picture in which {\em bubbles} replace corollas, {\em marked points} replace leaves, and the root is denoted as a point labeled $\infty$ (Figure~\ref{btp}b). Using the above notation, the composition $s \circ_i t$ is defined by fusing the $\infty$ of the bubble $s$ with the $i^{\rm th}$ marked point of $t$.  The branches of the tree are now identified with {\em double points}, the places where bubbles intersect.
\end{exmp}


\subsection{}
\begin{figure} [h]
\centering {\includegraphics {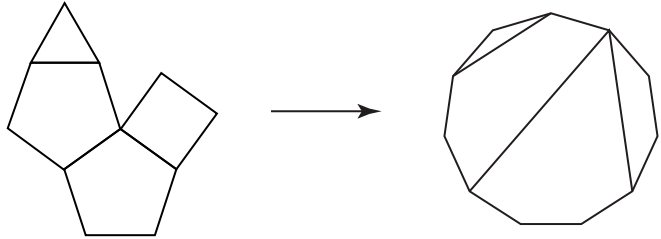}}
\caption{Polygon composition}
\label{onepoly}
\end{figure}

Taking yet another dual, we can define an operad structure on a collection of {\em polygons} (modulo an appropriate equivalence relation) as shown in Figure~\ref{btp}c.  Each bubble corresponds to a polygon, where the number of marked and double points become the number of sides; the fusing of points is associated with the gluing of faces.  The nicest feature of polygons is that, unlike corollas and bubbles, the iterated composition of polygons yields a polygon with marked diagonals (Figure~\ref{onepoly}).

\begin{figure} [h]
\centering {\includegraphics {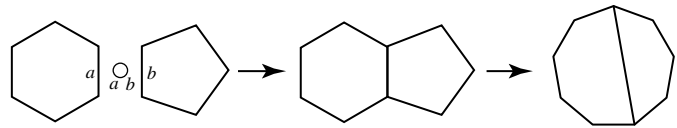}}
\caption{{\em Mosaic} composition}
\label{polycomp}
\end{figure}

Unlike the {\em rooted} trees, this {\em mosaic} operad is {\em cyclic} in the sense of Getzler and Kapranov~\cite[\S2]{cyclic}.  The most basic case (Figure~\ref{polycomp}) shows how two polygons, with sides labeled $a$ and $b$ respectively, compose to form a new polygon. The details of this operad are made precise in ~\S\ref{mosaic}.

\subsection{} \label{ss:lcubes}
\begin{figure} [h]
\centering {\includegraphics {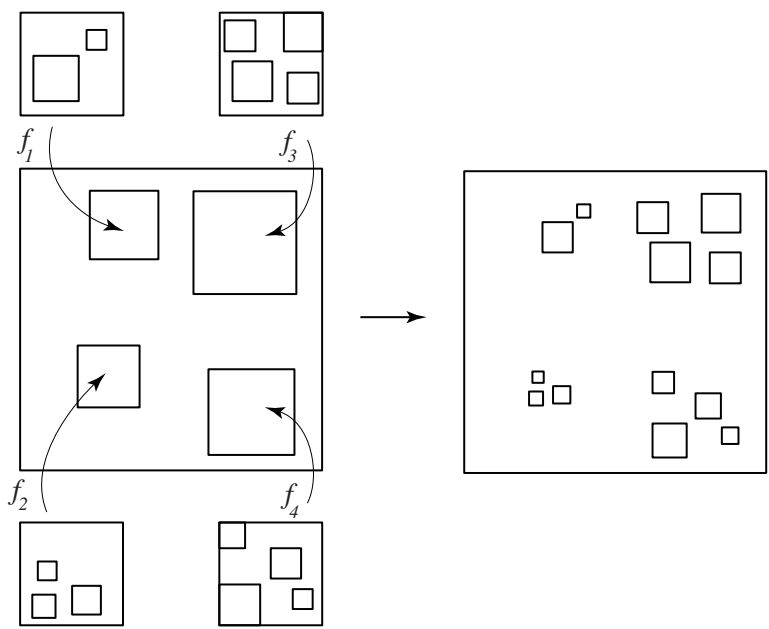}}
\caption{{\em Little cubes} composition}
\label{cubes}
\end{figure}

In the work of Boardman and Vogt~\cite[\S2.6]{bv}, an operad is presented using $m$ dimensional cubes $I^m \subset \R^m$.  An element $\Lc{n}$ of this {\em little cubes} operad is the space of an ordered collection of $n$ cubes linearly embedded by $f_i:I^m \hookrightarrow I^m$, with disjoint interiors and axes parallel to $I^m$.  The $f_i$'s are uniquely determined by the $2n$-tuple of points $(a_1, b_1, \ldots ,a_n, b_n)$ in $I^m$, corresponding to the images of the lower and upper vertices of $I^m$.  An element $\sg \in \Sg_n$ acts on $\Lc{n}$ by permuting the labeling of each cube:
$$(a_1, b_1, \ldots ,a_n, b_n) \mapsto (a_{\sg(1)}, b_{\sg(1)}, \ldots ,a_{\sg(n)}, b_{\sg(n)}).$$
The composition operation \eqref{e:operad} is defined by taking $n$ spaces $\Lc{k_i}$ (each having $k_i$ embedded cubes) and embedding them as an ordered collection into $\Lc{n}$. Figure~\ref{cubes} shows an example for the two dimensional case when $n = 4$.

Boardman showed that the space of $n$ distinct cubes in $\R^m$ is homotopically equivalent to $\Con^n(\R^m)$, the configuration space on $n$ distinct labeled points in $\R^m$.\footnote{The equivariant version of this theorem is proved by May in~\cite[\S4]{may1}.} When $m = 2$, $\Con^n(\R^2)$ is homeomorphic to $\C^n - \Delta$, where $\Delta$ is the {\em thick} diagonal $\{(x_1, \ldots , x_n) \in \C^n \: | \: \exists \: i, j, \: i \neq j \:,\:  x_i = x_j\}$.  Since the action of $\Sg_n$ on $\C^n - \Delta$ is free, taking the quotient yields another space $(\C^n - \Delta) / \Sg_n$.  It is well-known that both these spaces are aspherical, having all higher homotopy groups vanish~\cite{cdavis}.  The following short exact sequence of fundamental groups results:
$$ \pi_1 (\C^n - \Delta) \rightarrowtail \pi_1 ((\C^n - \Delta) / \Sg_n) \twoheadrightarrow \Sg_n.$$
But $\pi_1$ of $\C^n - \Delta$ is simply $\Pb_n$, the pure braid group.  Similarly, $\pi_1$ of $\C^n - \Delta$ quotiented by all permutations of labelings is the braid group $\B_n$.  Therefore, the short exact sequence above takes on the more familiar form:
$$\Pb_n \rightarrowtail \B_n \twoheadrightarrow \Sg_n.$$
We will return to these ideas in~\S\ref{quasi}.

%
%

\section {The Moduli Space}

\subsection{} \label{ss:collide}
The moduli space of Riemann spheres with $n$ punctures,
$${\mathcal M}_0^{n}(\C) = \Con^n(\C \Pj^1)/\Pj \Gl_2(\C),$$
has been studied extensively~\cite{keel}. It has a Deligne-Mumford-Knudsen compactification ${{\overline{\mathcal M}}{^n_0(\C)}}$, a smooth variety of complex dimension $n-3$.  In fact, this variety is defined over the integers; we will look at the {\em real} points of this space.  These are the set of fixed points of ${{\overline{\mathcal M}}{^n_0(\C)}}$ under complex conjugation.

\begin{defn}
The moduli space \M{n} of configurations of $n$ smooth points on punctured stable real algebraic curves of genus zero is a compactification of the quotient
$((\R \Pj^1)^n - \Delta)/\Pj \Gl_2(\R),$
where $\Delta$ is the thick diagonal.
\end{defn}

\begin{rem}
This is an action of a non-compact group on a non-compact space.  Geometric invariant theory gives a natural compactification for this quotient, defined combinatorially in terms of bubble trees or algebraically as a moduli space of real algebraic curves of genus zero with $n$ points, which are stable in the sense that they have only finitely many automorphisms.
\end{rem}

A point of \oM{n} can be visualized as a bubble (that is, $\R\Pj^1$) with $n$ {\em distinct} labeled points.  For a particular labeling, the configuration space of such points gives us a fundamental domain of \oM{n}.  There are $n!$ possible labelings.  However, since there exists a copy of the dihedral group $D_n$ in $\Pj \Gl_2(\R)$, and since \oM{n} is defined as a quotient by $\Pj \Gl_2(\R)$, two labeled bubbles are identified by an action of $D_n$.  Therefore, there are $\frac{1}{2}(n-1)!$ copies of the fundamental domain that make up \oM{n}.  Since we remove the thick diagonal, these domains are open cells.

\begin{figure}[h]
\centering {\includegraphics {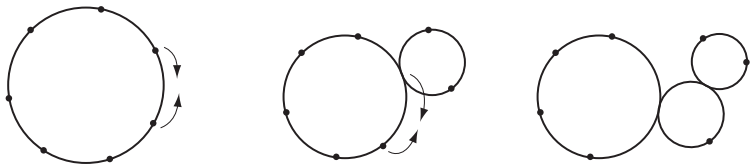}}
\caption{Collision on bubbles}
\label{bcollide}
\end{figure}

In \M{n}, however, these marked points are allowed to `collide' in the following sense: As two adjacent points $p_1$ and $p_2$ of the bubble come closer together and try to collide, the result is a new bubble fused to the old at the point of collision (a double point), where the marked points $p_1$ and $p_2$ are now on the new bubble (Figure~\ref{bcollide}).  Note that each bubble must have at least three marked or double points in order to be stable.

\begin{figure}[h]
\centering {\includegraphics {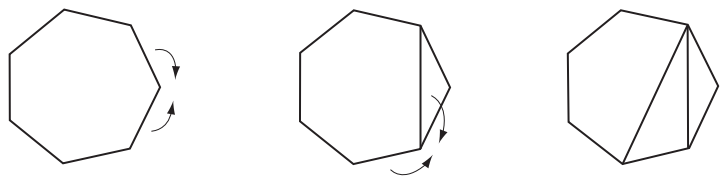}}
\caption{Collision on polygons}
\label{pcollide}
\end{figure}

The mosaic operad encapsulates all the information of the bubbles, enabling one to look at the situation above from the vantage point of polygons.  Having $n$ marked points on a circle now corresponds to an $n$-gon; when two adjacent sides $p_1$ and $p_2$ of the polygon try to collide, a diagonal of the polygon is formed such that $p_1$ and $p_2$ lie on one side of the diagonal (Figure~\ref{pcollide}). 

What is quite striking about \M{n} is that its homotopy properties are completely encapsulated in the fundamental group.

\begin{thm} \textup{\cite[\S5.1]{djs}}
\M{n} is aspherical.
\end{thm}

\noindent We will return to the structure of the fundamental group in~\S\ref{quasi}.


\subsection{} \label{mosaic}
We now turn to defining the mosaic operad and relating its properties with the structure of \oM{n}.  Let $S^1$ be the unit circle bounding $\D$, the disk endowed with the Poincar\'{e} metric; this orients the circle.  The geodesics in $\D$ correspond to open diameters of $S^1$ together with open circular arcs orthogonal to $S^1$.  The group of isometries on $\D$ is $\Pj\Gl_2(\R)$~\cite[\S4]{rat}.

A configuration of $n$ distinct points on $\R\Pj^1$ defines an {\em ideal} polygon in $\D$, with all vertices on the circle and geodesic sides. Let $\G(n,0)$ be the space of such configurations, modulo $\Pj\Gl_2(\R)$, and let $\G(n,k)$ be the space of such ideal polygons marked with $k$ non-intersecting geodesics between non-adjacent vertices. We want to think of the elements of $\G(n,k)$ as limits of configurations in $\G(n,0)$ in which $k$ sets of points have coalesced (see discussion above). Specifying $k$ diagonals defines a decomposition of an $n$-gon into $k+1$ smaller polygons, and we can topologize $\G(n,k)$ as a union of $(k+1)$-fold products of $\G(m,0)$'s corresponding to this decomposition.  For example, to the one dimensional space $\G(4,0)$ we attach zero dimensional spaces of the form $\G(3,0) \times \G(3,0)$. The combinatorics of these identifications can be quite complicated, but Stasheff's associahedra were invented to solve just such problems, as we will see in~\S\ref{ss:gl-mon} below.

Henceforth, we will visualize elements of $\G(n,k)$ as $n$-gons with $k$ non-intersecting diagonals, and we write $\G(n)$ for the space of $n$-gons with any number of such diagonals. Elements of $\G(n)$ inherit a natural cyclic order on their sides, and we write $\G^L(n)$ for the space of $n$-gons with labeled sides.


\begin{prop} \label{p:gln1}
There exists a bijection between the points of \,\oM{n} and the elements of \,$\G^L(n,0)$.
\end{prop}

\begin{rem}
Given an element in $\G(n,k)$, we can associate to it a dual tree.  Its vertices are the barycenters of polygons, defined using the Riemann mapping theorem, and the branches are geodesics between barycenters.  The leaves are geodesics that extend to points on $\R\Pj^1$ midway between two adjacent marked points on $\R\Pj^1$. It then follows that \M{n} is a space of {\em hyperbolic} planar trees. This perspective naturally gives a Riemann metric to \M{n}.
\end{rem}

\begin{defn}
Given $G \in \G^L(m,l)$ and $G_i \in \G^L(n_i,k_i)$ (where $1 \leq i \leq m$), there are composition maps
$$G \ _{a_{1}} \!\! \circ \! _{b_{1}} \  G_1 \ _{a_{2}} \!\! \circ \! _{b_{2}} \
\cdots \ _{a_{m}} \!\! \circ \! _{b_{m}} \  G_{m} \mapsto G_t,$$
where $G_t \in \G^L(-m + \sum n_i,\ m + l + \sum k_i)$.  The object $G_t$ is obtained by gluing side $a_i$ of $G$ along side $b_i$ of $G_i$.  The symmetric group $\Sg_n$ acts on $G_n$ by permuting the labeling of the sides.  These operations define the {\em mosaic} operad $\{\G^L(n,k)\}$.
\end{defn}

\begin{rem}
The one dimensional case of the little cubes operad is $\{\Li{n}\}$, the {\em little intervals} operad.  An element $\Li{n}$ is an ordered collection of $n$ embeddings of the interval $I \hookrightarrow I$, with disjoint interiors.  The notion of {\em trees} and {\em bubbles}, shown in Figure~\ref{btp}, is encapsulated in this intervals operad.  Furthermore, after embedding $I$ in $\R$ and identifying $\R \cup \infty$ with $\R\Pj^1$, the mosaic operad $\{\G^L(n,k)\}$ becomes a compactification of $\{\Li{n}\}$.
\end{rem}


\subsection{} \label{ss:gl-mon}
We now define the fundamental domain of \M{n} as a concrete geometric object and present its connections with the mosaic operad.

\begin{defn}
Let ${\mathcal A}$ be the space of $n-3$ distinct points $\{t_1, \ldots, t_{n-3}\}$ on the interval $[0,1]$ such that $0 < t_1 < \cdots < t_{n-3} <1$.  Identifying $\R \cup \infty$ with $\R\Pj^1$ carries the set $\{0, t_1, \ldots, t_{n-3}, 1, \infty\}$ of $n$ points onto $\R\Pj^1$. Therefore, there exists a natural inclusion of ${\mathcal A}$ in \M{n}.  Define the {\em associahedron} $K_{n-1}$ as the closure of the space ${\mathcal A}$ in \M{n}.
\end{defn}

\begin{prop} \label{p:gln2}
An interior point of \,$K_{n-1}$ corresponds to an element of \,$\G(n,0)$, and an interior point of a codim $k$ face corresponds to an element of \,$\G(n,k)$.
\end{prop}

\begin{proof}
Since $\Sg_3 \subset \Pj\Gl_2(\R)$, one can fix three of the $n$ distinct points on $\R\Pj^1$ to be $0, 1,$ and $\infty$.  Thus, the associahedron $K_{n-1}$ can be identified with the cell tiling $\M{n}$ and the proposition follows from the construction of $\G(n,k)$.
\end{proof}

The relation between the $n$-gon and $K_{n-1}$ is further highlighted by a work of Lee~\cite{lee}, where he constructs a polytope $Q_n$ that is dual to $K_{n-1}$, with one vertex for each diagonal and one facet for each triangulation of an $n$-gon.  He then proves the symmetry group of $Q_n$ to be the dihedral group $D_n$.  Restated, it becomes

\begin{prop} \textup{\cite[\S5]{lee}}
$D_n$ acts as a group of isometries on $K_{n-1}$.
\end{prop}

\begin{figure} [h]
\centering {\includegraphics {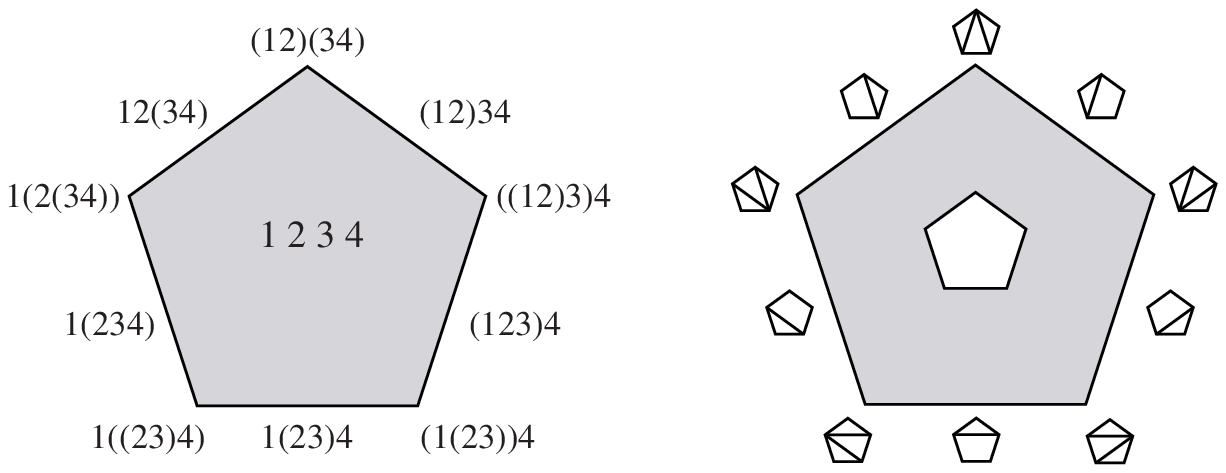}}
\caption{$K_4$}
\label{k4}
\end{figure}

\begin{hnote}
Stasheff classically defined the associahedron $K_{n-1}$ for use in homotopy theory~\cite[\S6]{jds} as a CW-ball with codim $k$ faces corresponding to using $k$ sets of parentheses meaningfully on ${n-1}$ letters.\footnote{From the definition above, the ${n-1}$ letters can be viewed as the points $\{0, t_1, \ldots, t_{n-3}, 1\}$.}  It is easy to describe the associahedra in low dimensions: $K_2$ is a point, $K_3$ a line, and $K_4$ a pentagon. The two descriptions of the associahedron, using polygons and parentheses, are compatible: Figure~\ref{k4} illustrates $K_4$ as an example. The associahedra have continued to appear in a vast number of mathematical fields, gradually acquiring more and more structure, {\em cf}.\ \cite{zie}.
\end{hnote}

\subsection{}
The polygon relation to the associahedron enables the use of the mosaic operad structure on $K_{n-1}$.

\begin{prop} \label{p:asdecomp} \textup{\cite[\S2]{jds}}  
Each face of $K_{n-1}$ is a product of lower dimensional associahedra.
\end{prop}

\noindent In general, the codim $k-1$ face of the associahedron $K_{m-1}$ will decompose as $$K_{n_1-1} \times \cdots \times K_{n_k-1} \hookrightarrow K_{m-1},$$
where $\sum n_i = m + 2(k-1)$ and $n_i \geq 3$.  This parallels the mosaic operad structure
$$G(n_1) \circ \cdots \circ G(n_k) \mapsto G(m),$$
where $G(n_i) \in \G^L(n_i,0),\ G(m) \in \G^L(m,k-1)$, and the gluing of sides is arbitrary.  Therefore, the product in Proposition~\ref{p:asdecomp} is indexed by the internal vertices of the tree corresponding to the face of the associahedron.

\begin{figure} [h]
\centering {\includegraphics {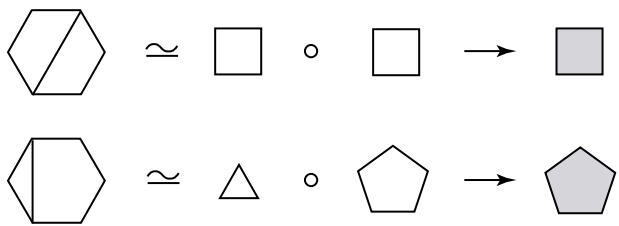}}
\caption{Codim one cells of $K_5$}
\label{k5codim1}
\end{figure}

\begin{exmp}
We look at the codim one faces of $K_5$. The three dimensional $K_5$ corresponds to a 6-gon, which has two distinct ways of adding a diagonal.  One way, in Figure~\ref{k5codim1}a, will allow the 6-gon to decompose into a product of two 4-gons ($K_3$'s).  Since $K_3$ is a line, this codim one face yields a square. The other way, in Figure~\ref{k5codim1}b, decomposes the 6-gon into a 3-gon ($K_2$) and a 5-gon ($K_4$).  Taking the product of a point and a pentagon results in a pentagon.
\end{exmp}

\begin{figure} [h]
\centering {\includegraphics {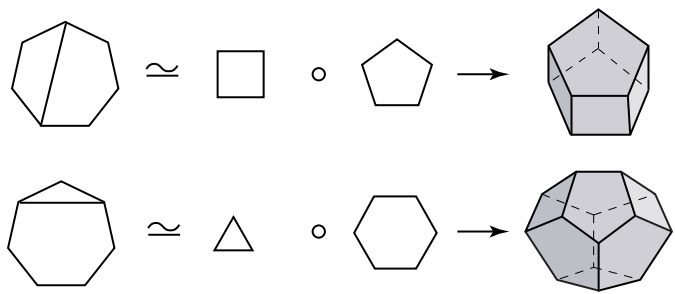}}
\caption{Codim one cells of $K_6$}
\label{k6codim1}
\end{figure}

\begin{exmp}
We look at the codim one faces of $K_6$.  Similarly, Figure~\ref{k6codim1} shows the decomposition of the codim one faces of $K_6$, a pentagonal prism and $K_5$.
\end{exmp}

%
%

\section {The Tessellation} \label{s:tess}

\subsection{} \label{twisting}
We extend the combinatorial structure of the associahedra to \M{n}.  Propositions~\ref{p:gln1} and \ref{p:gln2} show the correspondence between the associahedra in \M{n} and $\G^L(n,k)$.  We investigate how these copies of $K_{n-1}$ glue to form \M{n}.

\begin{figure} [h]
\centering {\includegraphics {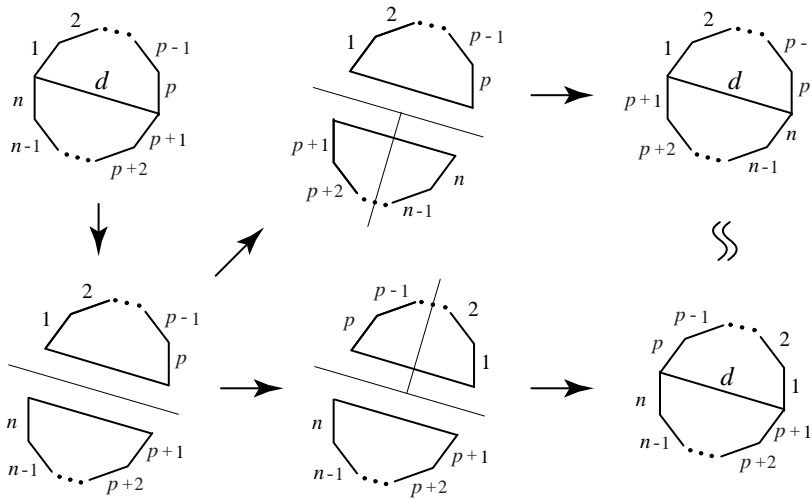}}
\caption{{\em Twist} along $d$}
\label{twist}
\end{figure}

\begin{defn}
Let $G \in \G^L(n,k)$ and $d$ be a diagonal of $G$.  A {\em twist} along $d$, denoted by $\nabla_d (G)$, is the element of $\G^L(n,k)$ obtained by `breaking' $G$ along $d$ into two parts, `twisting' one of the pieces, and `gluing' them back (Figure~\ref{twist}).
\end{defn}

\noindent The twisting operation is well-defined since the diagonals of an element in $\G^L(n,k)$ do not intersect. Furthermore, it does not matter which piece of the polygon is twisted since the two results are identified by an action of $D_n$. It immediately follows that $\nabla_d \cdot \nabla_d = e,$ the identity element.

\begin{prop} \label {p:twist}
Two elements, $G_1, G_2 \in \G^L(n,k)$, representing codim $k$ faces of associahedra, are identified in \M{n} if there exist diagonals $d_1, \ldots, d_r$ of $G_1$ such that $$(\nabla_{d_1} \cdots \nabla_{d_r}) (G_1) = G_2.$$
\end{prop}

\begin{figure} [h]
\centering {\includegraphics {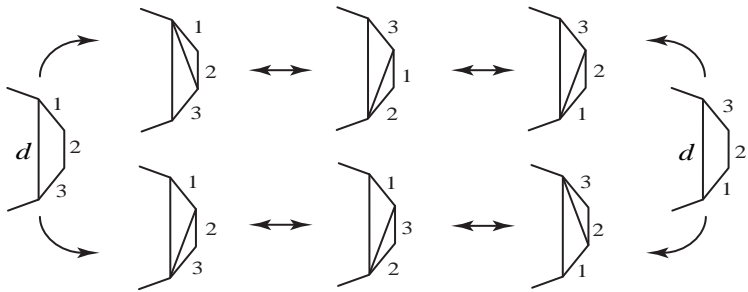}}
\caption{$\nabla_d$ does not affect the cell}
\label{twistpf}
\end{figure}

\begin{proof}
As two adjacent points $p_1$ and $p_2$ on $\R\Pj^1$ collide, the result is a new bubble fused to the old at a point of collision $p_3$, where $p_1$ and $p_2$ are on the new bubble. The location of the three points $p_i$ on the new bubble is {\em irrelevant} since $\Sg_3 \subset \Pj\Gl_2(\R)$.  In terms of polygons, this means $\nabla_d$ does not affect the cell, where $d$ is the diagonal representing the double point $p_3$.  In general, it follows that the labels of triangles can be permuted without affecting the cell.  Let $G$ be an $n$-gon with diagonal $d$ partitioning $G$ into a square and an $(n-2)$-gon.  Figure~\ref{twistpf} shows that since the square decomposes into triangles, the cell corresponding to $G$ is invariant under the action of $\nabla_d$.  Since any partition of $G$ by a diagonal $d$ can be decomposed into triangles, it follows by induction that $\nabla_d$ does not affect the cell.
\end{proof}

\begin{thm} \label{t:kxs}
There exists a surjection
$$K_{n-1} \times_{D_n} \Sg_n \rightarrow \M{n},$$
which is a bijection on the interior of the cells.
In particular, $\frac{1}{2}(n-1)!$ copies of $K_{n-1}$ tessellate \M{n}.
\end{thm}

\begin{proof}
The bijection on the interior of the cells follows immediately from the discussion in~\S\ref{ss:collide}. The map is not an injection since the boundaries of the associahedra are glued according to Proposition~\ref{p:twist}.
\end{proof}


\subsection{}
\begin{figure} [h]
\centering {\includegraphics {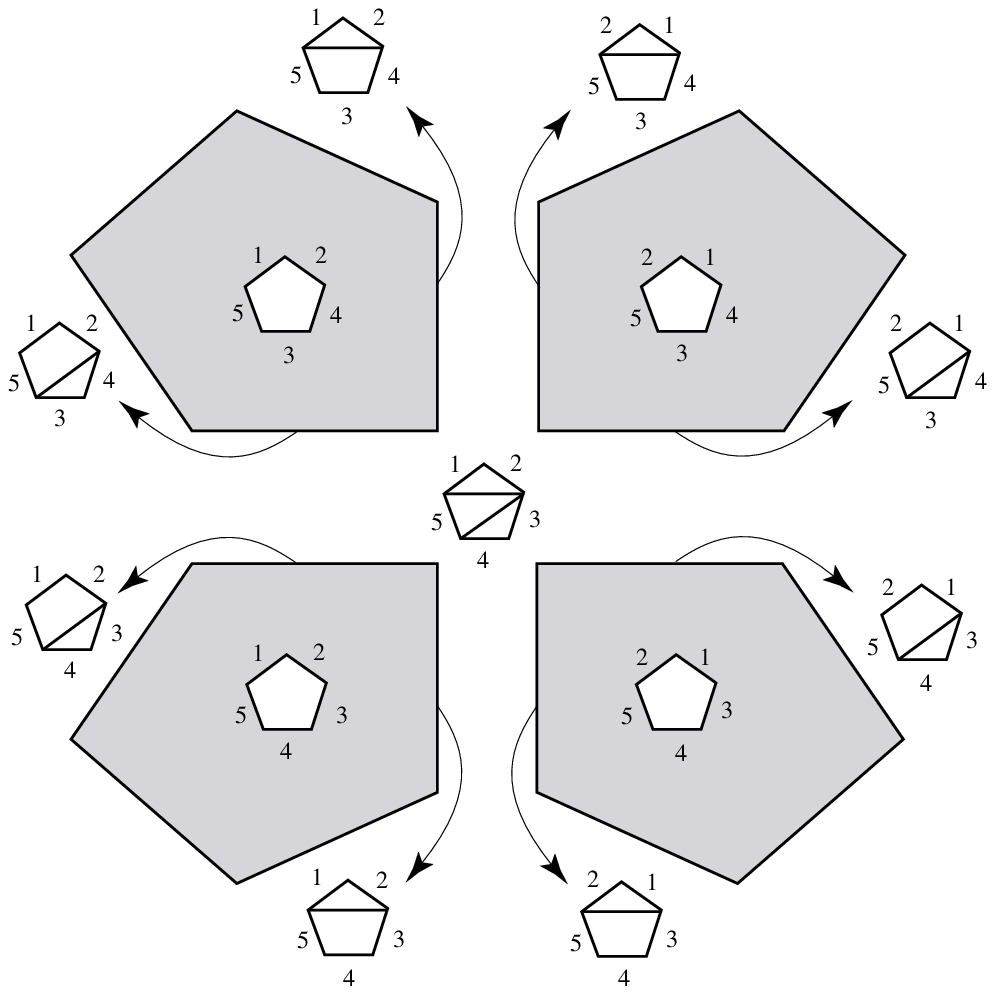}}
\caption{A piece of \M{5}}
\label{pieces}
\end{figure}

In Figure~\ref{pieces}, a piece of \M{5} represented by labeled polygons with diagonals is shown.  Note how two codim one pieces (lines) glue together and four codim two pieces (points) glue together.  Understanding this gluing now becomes a combinatorial problem related to $\G^L(n,k)$.

\begin{nota}
Let \Bnd{x}{\mathfrak X} be the number of codim $x$ cells in a CW-complex $\mathfrak X$.  For a fixed codim $y_2$ cell in \M{n}, and for $y_1 < y_2$, let \Cobnd{y_1}{y_2} be the number of codim $y_1$ cells in \M{n} whose boundary contains the codim $y_2$ cell.  Note the number \Cobnd{y_1}{y_2} is well-defined by Theorem~\ref{t:kxs}.
\end{nota}

\begin{lem} \label{l:cayley}
$$\Bnd{k}{K_{n-1}} = \frac{1}{k+1} \; \binom{n-3}{k} \; \binom{n-1+k}{k}.$$
\end{lem}

\begin{proof}
This is obtained by just counting the number of $n$-gons with $k$ non-intersecting diagonals, done by A. Cayley in 1891~\cite{cay}.
\end{proof}

\begin{lem} \label{l:codim}
$$\Cobnd{k-t}{k} = 2^t \; \binom{k}{t}.$$
\end{lem}

\begin{proof}
The boundary components of a cell corresponding to an element in $\G^L(n,k)$ are obtained by adding non-intersecting diagonals.  To look at the coboundary cells, diagonals need to be {\em removed}.  For each diagonal removed, two cells result (coming from the {\em twist} operation);  removing $t$ diagonals gives $2^t$ cells.  We then look at all possible ways of removing $t$ out of $k$ diagonals.
\end{proof}

\begin{thm} \label{t:euler}
\begin{equation}
\chi (\M{n}) = 
\begin{cases}
0 & n \text{ even}\\
(-1)^{\frac{n-3}{2}}(n-2)((n-4)!!)^2 & n \text{ odd.}
\end{cases}
\label{e:euler}
\end{equation}
\end{thm}

\begin{proof}
It is easy to show the following:
$$\Bnd{k}{\M{n}} \cdot \Cobnd{0}{k} = \Bnd{0}{\M{n}} \cdot \Bnd{k}{K_{n-1}}.$$
Using Theorem~\ref{t:kxs} and Lemmas~\ref{l:cayley} and~\ref{l:codim}, we solve for \Bnd{k}{\M{n}}; but this is simply the number of codim $k$ cells in \M{n}.  Therefore,
$$\chi (\M{n}) = \sum_{k=0}^{n-3} (-1)^{n-3-k} \;\; \frac{(n-1)!}{ 2^{k+1}} \:\; \frac{1}{k+1} \; \binom{n-3}{k} \; \binom{n-1+k}{k}.$$
This equation can be reduced to the desired form.
\end{proof}

\begin{rem}
Professor F.\ Hirzebruch has kindly informed us that he has shown, using techniques of Kontsevich and Manin~\cite{km}, that the signature of ${{\overline{\mathcal M}}{^n_0(\C)}}$ is given by \eqref{e:euler}. He remarks that the equivalence of this signature with the Euler number of the space of real points is an elementary consequence of the Atiyah-Singer $G$-signature theorem. 
\end{rem}

%
%

\section {The Hyperplanes}

\subsection{} \label{braidarr}
Another approach to \M{n} is from a {\em top-down} perspective using hyperplane arrangements as formulated by Kapranov~\cite[\S4.3]{kapchow} and described by Davis, Januszkiewicz, and Scott~\cite[\S0.1]{djs}.

\begin{defn}
Let $V^n \subset \R^{n-1}$ be the hyperplane defined by $\Sigma x_i = 0$.  For $1 \leq i < j \leq n-1$, let $H^n_{ij} \subset V^n$ be the hyperplane defined by $x_i = x_j$. The {\em braid arrangement} is the collection of subspaces of $V^n$ generated by all possible intersections of the $H^n_{ij}$.
\end{defn}

If $\Hu^n$ denotes the collection of subspaces $\{H^n_{ij}\}$, then $\Hu^n$ cuts $V^n$ into $(n-1)!$ simplicial cones. Let $\Sg(V^n)$ be the sphere in $V^n$ and let $\Pj(V^n)$ be the projective sphere in $V^n$ (that is, $\R\Pj^{n-3}$).  Let \Ba{n} to be the intersection of $\Hu^n$ with $\Pj(V^n)$; the arrangement \Ba{n} cuts $\Pj(V^n)$ into $\frac{1}{2}(n-1)!$ open $n-3$ simplices.  

\begin{defn}
Let $\ba^k$ be a codim $k$ {\em irreducible} cell of $\Pj(V^n)$ if $\binom{k+1}{2}$ hyperplanes of $\Hu^n$ intersect there.\footnote{The use of the word {\em irreducible} comes from \cite{djs} in reference to Coxeter groups.}
\end{defn}

\begin{figure} [h]
\centering {\includegraphics {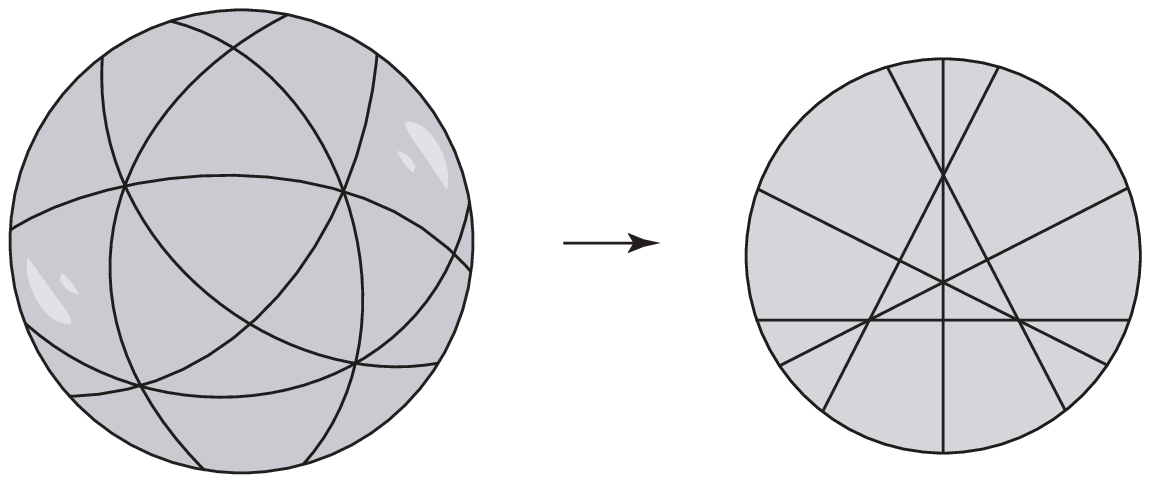}}
\caption{\protect{$\Sg(V^5) \rightarrow \Pj(V^5)$}}
\label{svpv}
\end{figure}

\begin{exmp}
We look at the case when $n=5$. Figure~\ref{svpv} shows the `scars' on the manifolds made by $\Hu^5$.  On $\Pj(V^5)$, there are four places where three hyperplanes intersect, corresponding to the four codim two irreducible points.
\end{exmp}

\begin{defn}
Replace $\ba^k$ with $\sba^k$, the sphere bundle associated to the normal bundle of $\ba^k \subset \Pj(V^n)$.  This process yields a manifold with boundary.  Then projectify $\sba^k$ into $\pba^k$, the projective sphere bundle.  This defines a manifold without boundary, called the {\em blow-up of \,$\Pj(V^n)$ along $\ba^k$}.
\end{defn}  

\begin{rem}
Replacing $\ba^k$ with $\sba^k$ for {\em any} dimension $k$ creates a {\em new} manifold with boundary.  However, blowing up along $\ba^{k}$ defines a new manifold for all dimensions {\em except} codim one.  That is, for codim one, projectifying $\sba^k$ into $\pba^k$ annuls the process of replacing $\ba^k$ with $\sba^k$.
\end{rem}

\begin{prop} \textup{\cite[\S4.3]{kapchow}} \label{p:kap}
The iterated blow-up of \,$\Pj(V^n)$ along the cells $\{\ba^k\}$ in {\em increasing} order of dimension yields \M{n}. It is inessential to specify the order in which cells $\{\ba^k\}$ of the {\em same} dimension are blown up.
\end{prop}

Therefore, the compactification of \oM{n} is obtained by replacing the set $\{\ba^k\}$ with $\{\pba^k\}$.  The {\em closure} of \oM{n} in $\Pj(V^n)$ is obtained by replacing the set $\{\ba^k\}$ with \{$\sba^k$\}; this procedure truncates each $n-3$ simplex of $\Pj(V^n)$ into the associahedron $K_{n-1}$.  We explore this method of truncation in~\S\ref{ss:truncate}.

\begin{figure} [h]
\centering {\includegraphics {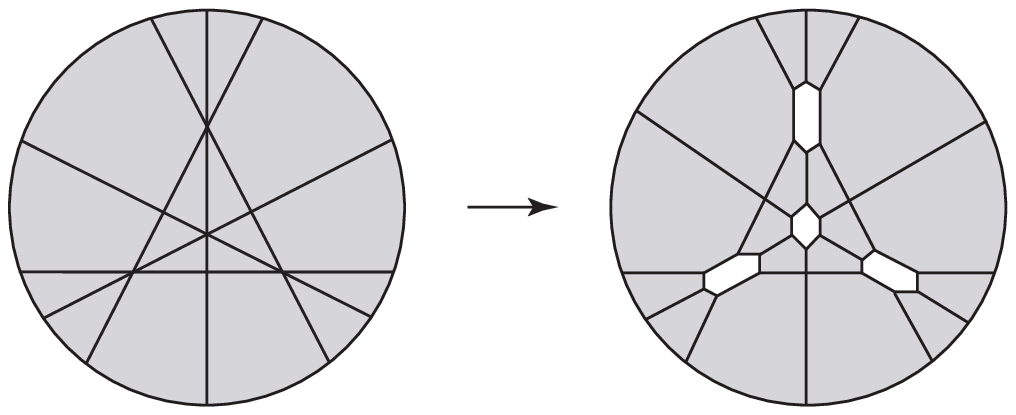}}
\caption{\protect{$\Pj(V^5) \rightarrow \M{5}$}}
\label{pvm05}
\end{figure}

\begin{exmp} \label {e:m05blowup}
The blow-up of $\Pj(V^5)$ yielding \M{5} is shown in Figure~\ref{pvm05}. The arrangement \Ba{5} on $\Pj(V^5) \simeq \R\Pj^2$ yields six lines forming twelve $2$-simplices;  the irreducible components of codim two turn out to be the points $\{\ba^2_1, \ldots, \ba^2_4\}$ of triple intersection.  Blowing up along these components, we get $S^1$ as a hexagon for $\sba^2_i$ and $\R\Pj^1$ as a triangle for $\pba^2_i$. The associahedron $K_4$ is a pentagon, and the space \M{5} becomes tessellated by twelve such cells (shaded), an ``evil twin'' of the dodecahedron.  \M{5} appears as the connected sum of five real projective planes.
\end{exmp}

\begin{hnote}
The diagram of \M{5} shown in Figure~\ref{pvm05} is first found in a different context by Brahana and Coble in $1926$~\cite[\S1]{bc} relating to possibilities of maps with twelve five-sided countries.  
\end{hnote}


\subsection{}
Another way of looking at the moduli space comes from observing the inclusion $\Sg_3 \subset \Pj \Gl_2 (\R)$.  Since \M{n} is defined as $n$ distinct points on $\R \Pj^1$ quotiented by $\Pj\Gl_2 (\R)$, one can fix three of these points to be $0, 1,$ and $\infty$.  From this perspective we see that \M{3} is a point.  When $n=4$, the {\em cross-ratio} is a homeomorphism from \M{4} to $\R\Pj^1$,  the result of identifying three of the four points with $0, 1,$ and $\infty$.  In general, \M{n} becomes a manifold blown up from an $n-3$ dimensional torus, coming from the $(n-3)$-fold products of $\R \Pj^1$. Therefore, the moduli space {\em before} compactification can be defined as 
$$((\R \Pj^1)^n - \Delta^*)/\Pj \Gl_2(\R),$$
where $\Delta^* = \{(x_1, \ldots , x_n) \in (\R \Pj^1)^n \:|\: $at least 3 points collide\}.
Compactification is accomplished by blowing up along $\Delta^*$.

\begin{figure} [h]
\centering {\includegraphics {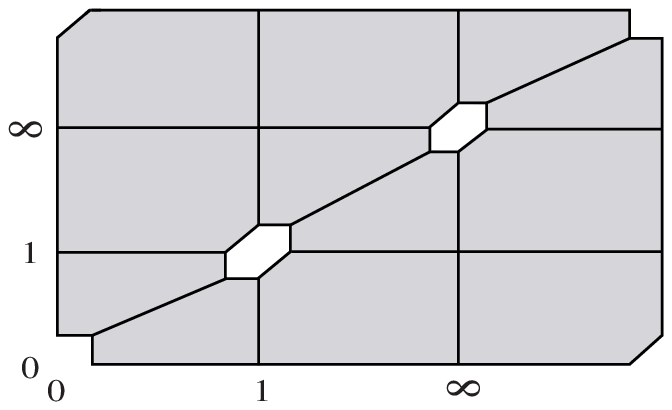}}
\caption{\M{5} from the torus} 
\label{m05c}
\end{figure}

\begin{exmp}
An illustration of \M{5} from this perspective appears in Figure~\ref{m05c}.  From the five marked points on $\R \Pj^1$, three are fixed leaving two dimensions to vary, say $x_1$ and $x_2$.  The set $\Delta$ is made up of seven lines $\{x_1, x_2 = 0, 1, \infty\}$ and $\{x_1 = x_2\}$, giving a space tessellated by six squares and six triangles.  Furthermore, $\Delta^*$ becomes the set of three points $\{x_1=x_2 = 0,1,\infty\}$; blowing up along these points yields the space \M{5}  tessellated by twelve pentagons.  This shows \M{5} as the connected sum of a torus with three real projective planes.
\end{exmp}

\begin{figure} [h]
\centering {\includegraphics {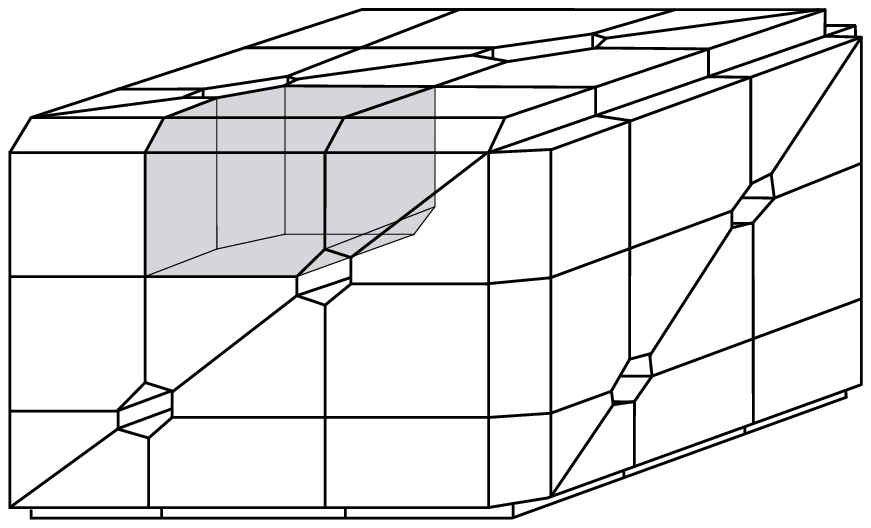}}
\caption{\M{6}}
\label{m06c}
\end{figure}

\begin{exmp} \label{e:m06}
In Figure~\ref{m06c}, a rough sketch of \M{6} is shown as the blow-up of a three torus. The set $\Delta^*$ associated to \M{6} has ten lines \{$x_i=x_j=0,1,\infty$\} and \{$x_1=x_2=x_3$\}, and three points \{$x_1=x_2=x_3=0,1,\infty$\}.  The lines correspond to the hexagonal prisms, nine cutting through the faces, and the tenth (hidden) running through the torus from the bottom left to the top right corner.  The three points correspond to places where four of the prisms intersect.

The shaded region has three squares and six pentagons as its codim one faces.  In fact, all the top dimensional cells that form \M{6} turn out to have this property; these cells are the associahedra $K_5$ (see Figure~\ref{k6codim1}b).
\end{exmp}

\subsection{}
We now introduce a construction which clarifies the structure of \M{n}.

\begin{defn} \textup{\cite[\S4]{kap}}
A double cover of \M{n}, denoted by \dM{n}, is obtained by fixing the $n^{\rm th}$ marked point on $\R\Pj^1$ to be $\infty$ and assigning it an orientation.\footnote{Kapranov uses the notation $\tilde S^{n-3}$ to represent this double cover.}
\end{defn}

\begin{figure} [h]
\centering {\includegraphics {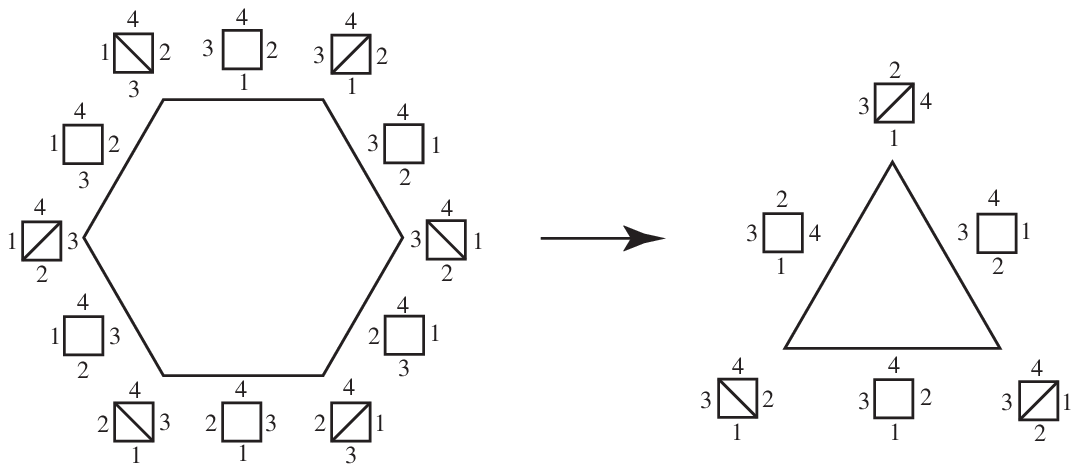}}
\caption{\protect{$\dM{4} \rightarrow \M{4}$}}
\label{m04}
\end{figure}

\begin{exmp}
Figure~\ref{m04} shows the polygon labelings of \dM{4} and \M{4}, being tiled by six and three copies of $K_3$ respectively.  In this figure, the label $4$ has been set to $\infty$.  Note that the map $\dM{4} \rightarrow \M{4}$ is the antipodal quotient.
\end{exmp}

The double cover can be constructed using blow-ups similar to the method described above; instead of blowing up the projective sphere $\Pj(V^n)$, we blow-up the sphere $\Sg(V^n)$. Except for the anomalous case of \dM{4}, the double cover is a {\em non-orientable} manifold. Note also that the covering map $\dM{n} \rightarrow \M{n}$ is the antipodal quotient, coming from the map $\Sg(V^n) \rightarrow \Pj(V^n)$.  Being a double cover, \dM{n} will be tiled by $(n-1)!$ copies of $K_{n-1}$.\footnote{These copies of $K_{n-1}$ are in bijection with the vertices of the {\em permutohedron} $P_{n-1}$~\cite{kap}.}  It is natural to ask how these copies glue to form \dM{n}.

\begin{defn}
A {\em marked twist} of an $n$-gon $G$ along its diagonal $d$, denoted by $\widetilde \nabla_d (G)$, is the polygon obtained by breaking $G$ along $d$ into two parts, reflecting the piece that does {\em not} contain the side labeled $\infty$, and gluing them back together. 
\end{defn}

The two polygons at the right of Figure~\ref{twist} turn out to be {\em different} elements in \dM{n}, whereas they are identified in \M{n} by an action of $D_n$.  The following is an immediate consequence of the above definitions and Theorem~\ref{t:kxs}.

\begin{cor} \label{c:kxs}
There exists a surjection
$$K_{n-1} \times_{\Z_n} \Sg_n \rightarrow \dM{n}$$
which is a bijection on the interior of the cells.
\end{cor}

\begin{rem}
The spaces on the left define the classical $A_{\infty}$ operad~\cite[\S2.9]{cyclic}.
\end{rem}

\begin{thm}
The following diagram is commutative:
$$\begin{CD}
(K_{n-1} \times \Sg_n)/_{\widetilde \nabla} @>>> \dM{n}\\
@VVV     @VVV\\
(K_{n-1} \times \Sg_n)/_{\nabla} @>>> \M{n}
\end{CD}$$
where the vertical maps are antipodal identifications and the horizontal maps are a quotient by $\Z_n$.
\end{thm}

\begin{proof}
Look at $K_{n-1} \times \Sg_n$ by associating to each $K_{n-1}$ a particular labeling of an $n$-gon.  We obtain $(K_{n-1} \times \Sg_n)/_{\widetilde \nabla}$ by gluing the associahedra along codim one faces using $\widetilde \nabla$ (keeping the side labeled $\infty$ fixed).  It follows that two associahedra will {\em never} glue if their corresponding $n$-gons have $\infty$ labeled on different sides of the polygon.  This partitions $\Sg_n$ into $\Sg_{n-1} \cdot \Z_n$, with each element of $\Z_n$ corresponding to $\infty$ labeled on a particular side of the $n$-gon.  Furthermore, Corollary~\ref{c:kxs} tells us that each set of the $(n-1)!$ copies of $K_{n-1}$ glue to form \dM{n}.  Therefore,
$(K_{n-1} \times \Sg_n)/_{\widetilde \nabla} \:=\: (K_{n-1} \times \Sg_{n-1})/_{\widetilde \nabla} \times \Z_n \:=\: \dM{n} \times \Z_n.$
\end{proof}

%
%

\section{The Blow-Ups}

\subsection{} \label{ss:observe}
The spaces \M{n} and $\R\Pj^{n-3}$ differ only by blow-ups, making the study of their structures crucial. Looking at the arrangement \Ba{n} on $\Pj(V^n)$, there turn out to be $n-1$ irreducible points $\{\ba^{n-3}\}$ in {\em general position}.  In other words, these points can be thought of as vertices of an $n-3$ simplex with an additional point at the center. Between every two $\ba^{n-3}$ points of \Ba{n}, there exists a $\ba^{n-4}$ line, resulting in $\binom{n-1}{n-3}$ such irreducible lines.  In general, $k$ irreducible points of \Ba{n} span a \mbox{$k-1$} dimensional irreducible cell; restating this, we get

\begin{prop} \label{p:icells}
The number of irreducible components $\ba^k$ in \Ba{n} equals
\begin{equation}
\binom{n-1}{k+1}.
\label{e:countirr}
\end{equation}
\end{prop}

\noindent The construction of the braid arrangement shows that around a point $\ba^{n-3}$ of $\Pj(V^n)$, the structure of \Ba{n} resembles the barycentric subdivision of an $n-3$ simplex.  We look at some concrete examples to demonstrate this.  

\begin{exmp}
In the case of \M{5}, Figure~\ref{pvm05}a shows the $\ba^2$ cells in general position; there are four points, three belonging to vertices of a $2$-simplex, and one in the center of this simplex.   Between every two of these points, there exists a $\ba^1$; we see six such lines.  Since these lines are of codim one, they need not be blown up.
Figure~\ref{pvm05}b shows the structure of a blown up point $\ba^2$ in \M{5}.  Notice that $\sba^2$ is a hexagon and $\pba^2$ is a triangle.  It is no coincidence that these correspond exactly to \dM{4} and \M{4} (see Figure~\ref{m04}).
\end{exmp}

\begin{exmp}
For the three dimensional \M{6}, the $\ba^3$ cells {\em and} the $\ba^2$ cells need to be blown up, {\em in that order}. Choose a codim three cell $\ba^3$; a neighborhood around $\ba^3$ will resemble the barycentric subdivision of a $3$-simplex.  Figure~\ref{braid6} shows four tetrahedra, each being made up of six tetrahedra (some shaded), pulled apart in space such that when glued together the result will constitute the aforementioned subdivision.  The barycenter is the point $\ba^3$.

\begin{figure} [h]
\centering {\includegraphics {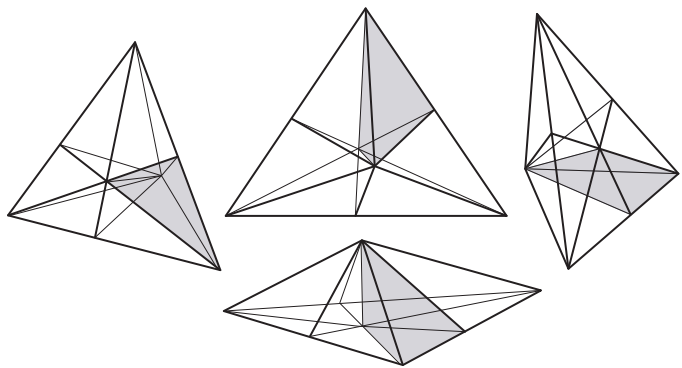}}
\caption{Barycentric subdivision of a $3$-simplex}
\label{braid6}
\end{figure}

The left-most piece of Figure~\ref{blow6} shows one of the tetrahedra from Figure~\ref{braid6}.  The map $f_1$ takes the barycenter $\ba^3$ to $\sba^3$ whereas the map $f_2$ takes each $\ba^2$ going through the barycenter to $\sba^2$.  When looking down at the resulting `blown up' tetrahedron piece, there are six pentagons (shaded) with a hexagon hollowed out in the center.  Taking $\sba^2$ to $\pba^2$ turns these hexagons into triangles.
\begin{figure} [h]
\centering {\includegraphics {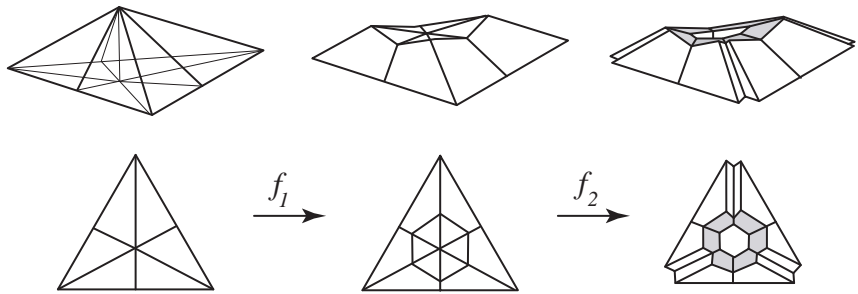}}
\caption{Blow-up of vertex and lines}
\label{blow6}
\end{figure}
Putting the four `blown up' tetrahedra pieces together, the faces of $\sba^3$ make up a two dimensional sphere tiled by 24 pentagons, with 8 hexagons (with antipodal maps) cut out. This turns out to be \dM{5}; projectifying $\sba^3$ to $\pba^3$ yields \M{5} as shown in Figure~\ref{dm05}.

\begin{figure} [h]
\centering {\includegraphics {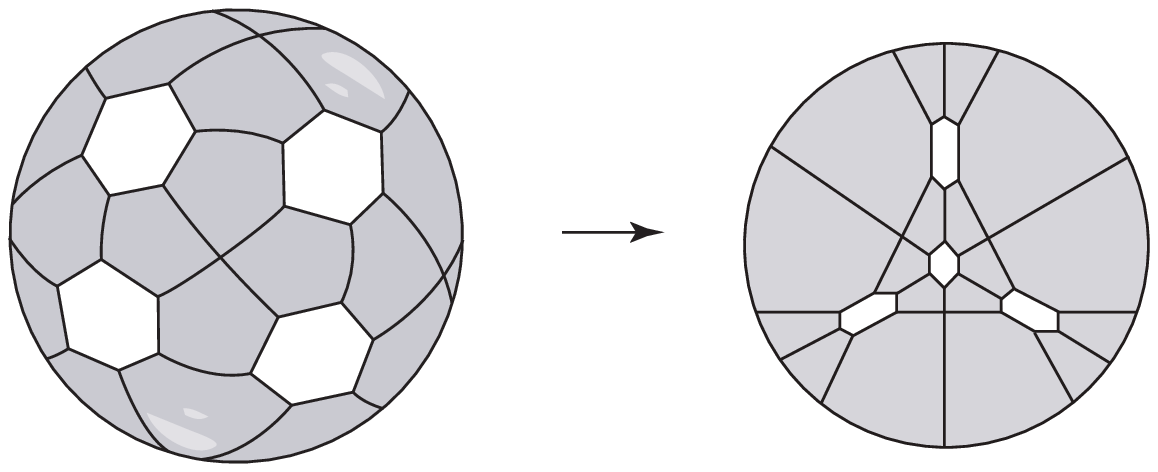}}
\caption{\protect{$\dM{5} \rightarrow \M{5}$}}
\label{dm05}
\end{figure}
\end{exmp}

This pattern seems to indicate that for \M{n}, blowing up along the point $\ba^{n-3}$ will yield \M{n-1}.  But what happens in general, when a codim $k$ cell $\ba^k$ is blown up? A glimpse of the answer was seen above with regard to the hexagons and triangles showing up in \M{6}.


\subsection{}
To better understand \M{n}, we analyze the structure of $\ba^k \in \Pj(V^n)$ before blow-ups and $\pba^k \in \M{n}$ after blow-ups.  This is done through the eyes of mosaics, looking at the faces of associahedra surrounding each blown up component of \Ba{n}.  The following is a corollary of Proposition~\ref{p:icells}.

\begin{cor} \label{c:icells}
Each irreducible cell \,$\ba^k$ corresponds to a choice of $k+1$ elements from the set $\{1, \ldots, n-1\}$.
\end{cor}

Choose an arbitrary $\ba^k$ and assign it such a choice, say $\{p_1, \ldots, p_{k+1}\}$, where $p_i \in \{1, \ldots, n-1\}$. We can think of this as an $n$-gon having a diagonal $d$ partitioning it such that $k+1$ labeled sides $\{p_1, \ldots, p_{k+1}\}$ lie on one side and $n-k-1$ labeled sides $\{p_{k+2}, \ldots, p_{n-1}, n\}$ lie on the other. Using the mosaic operad structure, $d$ decomposes the $n$-gon into $G_1 \circ \,G_2$, where $G_1 \in \G^L(k+2)$ and $G_2 \in \G^L(n-k)$, with the new sides $d_i$ of $G_i$ coming from $d$.  Note that $G_1 \circ \,G_2$ corresponds to the product of associahedra $K_{k+1} \times K_{n-k-1}$.

There are $(k+1)!$ different ways in which $\{p_1, \ldots, p_{k+1}\}$ can be arranged to label $G_1$.  However, since {\em twisting} is allowed along $d_1$, we get $\frac{1}{2}(k+1)!$ different labelings of $G_1$, each corresponding to a $K_{k+1}$.  But observe that this is {\em exactly} how one gets \M{k+2}, where the associahedra glue as defined in \S\ref{twisting}.  Therefore, a fixed labeling of $G_2$ gives $\M{k+2} \times K_{n-k-1}$; all possible labelings result in

\begin{thm}
In \,\M{n}, each irreducible cell \,$\ba^k$ in \,\Ba{n} becomes
\begin{equation}
\M{k+2} \times \M{n-k}.
\label{e:mxm}
\end{equation}
\end{thm}

\begin{exmp}
Since \M{3} is a point, the blown up $\ba^{n-3}$ cell becomes \M{n-1}, matching the earlier observations of \S\ref{ss:observe}.  Furthermore, \eqref{e:countirr} shows there to be \mbox{$n-1$} such structures.
\end{exmp}

\begin{exmp}
Although blowing up along codim one components does not affect the resulting manifold, we observe their presence in \M{5}.  From \eqref{e:countirr}, we get six such $\ba^1$ cells which become \M{4} after blow-ups.  The \M{4}'s are seen in Figure~\ref{pvm05} as the six lines cutting through $\R\Pj^2$.  Note that every line is broken into six parts, each part being a $K_3$.
\end{exmp}

\begin{exmp}
The space \M{6}, illustrated in Figure~\ref{m06c}, moves a dimension higher.\footnote{Although this figure is not constructed from the braid arrangement, it is homeomorphic to the structure described by the braid arrangement.}  There are ten $\ba^2$ cells, each becoming $\M{4} \times \M{4}$. These are the hexagonal prisms that cut through the three torus as described in Example~\ref{e:m06}.
\end{exmp}


\subsection{}
The question arises as to {\em why} $\, \M{n-k}$ appears in \M{n}.  The answer lies in the braid arrangement of hyperplanes.  Taking \M{6} as an example, blowing up along each point $\ba^3$ in \Ba{6} uses the following procedure:  A small spherical neighborhood is drawn around $\ba^3$ and the inside of the sphere is removed, resulting in $\sba^3$.  Observe that this sphere (which we denote as $\Sh$) is engraved with great arcs coming from \Ba{6}.  Projectifying, $\sba^3$ becomes $\pba^3$, and $\Sh$ becomes the projective sphere $\Pj\Sh$. Amazingly, the engraved arcs on $\Pj\Sh$ are \Ba{5}, and $\Pj\Sh$ can be thought of as $\Pj(V^5)$.  Furthermore, blowing up along the lines $\ba^2$ of \Ba{6} corresponds to blowing up along the points $\ba^2$ of \Ba{5} in $\Pj\Sh$.  As before, this new etching on $\Pj\Sh$ translates into an even lower dimensional braid arrangement, \Ba{4}.

It is not hard to see how this generalizes in the natural way: For \M{n}, the iterated blow-ups along the cells $\{\ba^{n-3}\}$ up to $\{\ba^2\}$ in turn create braid arrangements within braid arrangements.  Therefore, $\M{n-k}$ is seen in $\M{n}$.


\subsection{} \label{ss:truncate}
So far we have been looking at the structure of the irreducible cells $\ba^k$ before and after the blow-ups.  We now study how the $n-3$ simplex (tiling $\Pj(V^n)$) is truncated by blow-ups to form $K_{n-1}$ (tiling \M{n}).\footnote{For a detailed construction of this truncation from another perspective, see Appendix B of~\cite{jds2}.} Given a regular $n$-gon with one side marked $\infty$, define $\Gc$ to be the set of such polygons with one diagonal.

\begin{figure} [h]
\centering {\includegraphics {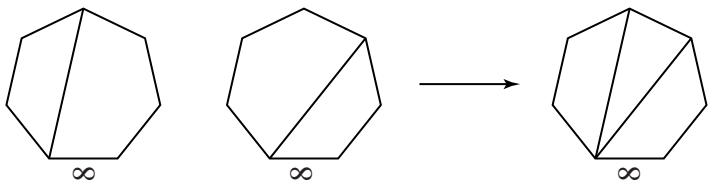}}
\caption{{\em Superimpose}}
\label{f:simpose}
\end{figure}

\begin{defn}
For $G_1, G_2 \in \Gc$, create a new polygon $G_{1,2}$ (with {\em two} diagonals) by {\em superimposing} the images of $G_1$ and $G_2$ on each other (Figure~\ref{f:simpose}).  $G_1$ and $G_2$ are said to satisfy the {\em \SI\ condition} if $G_{1,2}$ has non-intersecting diagonals.
\end{defn}

\begin{rem}
It follows from \S\ref{ss:gl-mon} that elements of $\Gc$ correspond bijectively to the codim one faces of $K_{n-1}$.  They are {\em adjacent} faces in $K_{n-1}$ if and only if they satisfy the \SI\ condition.  Furthermore, the codim two cell of intersection in $K_{n-1}$ corresponds to the superimposed polygon.  
\end{rem}

The diagonal of each element $G_i \in \Gc$ partitions the $n$-gon into two parts, with one part {\em not} having the $\infty$ label; call this the {\em free part of $G_i$}.  Define the set $\Gc^i$ to be elements of $\Gc$ having $i$ sides on their free parts. It is elementary to show that the order of $\Gc^i$ is $n-i$ (for $1 < i < n-1$).  In particular, the order of $\Gc^2$ is $n-2$, the number of sides (codim one faces) of an $n-3$ simplex.  Arbitrarily label each face of the simplex with an element of $\Gc^2$.

\begin{rem}
For some adjacent faces of the $n-3$ simplex, the \SI\ condition is not satisfied. This is an obstruction of the simplex in becoming $K_{n-1}$.  As we continue to truncate the cell, more faces will begin to satisfy the \SI\ condition.  We note that once a particular labeling is chosen, the labels of all the new faces coming from truncations (blow-ups) will be forced.  
\end{rem}

\begin{figure} [h]
\centering {\includegraphics {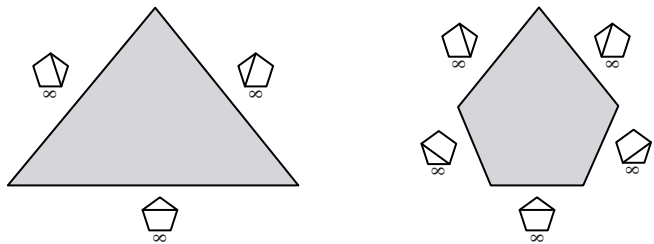}}
\caption{Truncation of $K_4$ by blow-ups}
\label{f:trunk4}
\end{figure}

When the zero dimensional cells are blown up, two vertices of the simplex are truncated. The labeling of the two new faces corresponds to the two elements of $\Gc^{n-2}$.  We choose the vertices and the labels such that the \SI\ condition is satisfied with respect to the {\em new} faces and {\em their} adjacent faces.  Figure~\ref{f:trunk4} shows the case for the $2$-simplex and $K_4$ (compare with Figure~\ref{pvm05}).

\begin{figure} [h]
\centering {\includegraphics {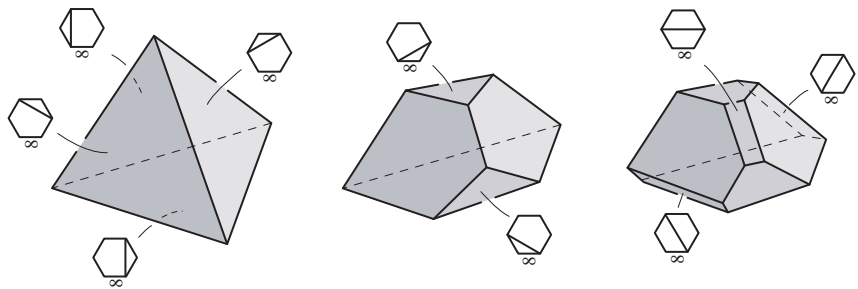}}
\caption{Truncation of $K_5$ by blow-ups}
\label{f:trunk5}
\end{figure}

The blow-up of one dimensional cells results in the truncation of three lines.  As before, the labels of the three new faces correspond to the three elements of $\Gc^{n-3}$, choosing edges and the labels such that the \SI\ condition is satisfied with respect to the new faces and their adjacent faces.  Figure~\ref{f:trunk5} shows the case for the $3$-simplex and $K_5$ (compare with Figures~\ref{braid6} and~\ref{blow6}). 

As we iterate the blow-ups in Proposition~\ref{p:kap}, we jointly truncate the \mbox{$n-3$} simplex using the above process.  The blow-ups of the codim $k$ irreducible cells $\ba^k$ add \mbox{$n-k-1$} new faces to the polytope, each labeled with an element from $\Gc^{k+1}$.  Note that Corollary~\ref{c:icells} is in agreement with this procedure:  Each irreducible cell $\ba^k$ corresponds to a choice of $k+1$ labels which are used on the elements of $\Gc^{k+1}$.
In the end, we are left with $\sum |\Gc^i|$ faces of the truncated polytope, matching the number of codim one faces of $K_{n-1}$.

%
%

\section{The Fundamental Group} \label{quasi}

\subsection{}
Coming full circle, we look at connections between the little cubes and the mosaic operads.  We would like to thank M.\ Davis, T.\ Januszkiewicz, and R.\ Scott for communicating some of their results in preliminary form~\cite{djs2}.  Their work is set up in the full generality of Coxeter groups and reflection hyperplane arrangements, but we explain how it fits into the notation of polygons and diagonals.

\begin{figure} [h]
\centering {\includegraphics {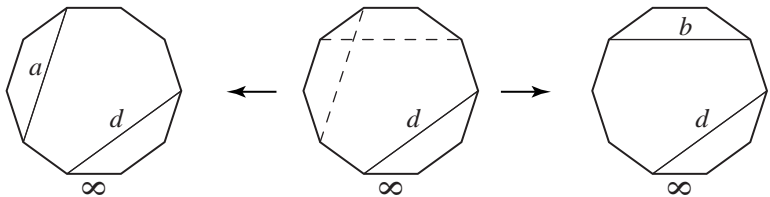}}
\caption{{\em Conjugate}}
\label{f:conjugate}
\end{figure}

\begin{defn}
Let $G_a, G_d \in \Gc$, with diagonals $a, d$ respectively, satisfy the \SI\ condition.  Let $G_b$ be the element in $\Gc$ after removing diagonal $d$ from $\widetilde \nabla_d (G_{a,d})$. We then say that $G_a$ and $G_b$ are {\em conjugate in $G_d$}.  Figure~\ref{f:conjugate} shows such a case.
\end{defn}

\begin{defn}
Let \Cox{n-1} be a group generated by elements $\{s_i\}$, in bijection with the elements $\{G_i\}$ of $\Gc$, with the following relations:

\vspace{3pt}
\begin{tabular}{cl} 
$s_i^2 = 1$ & \\
$s_d s_a = s_b s_d$ & if $G_a$ and $G_b$ are conjugate in $G_d$ \\
$s_a s_b = s_b s_a$ & if $G_a$ and $G_b$ satisfy the \SI\ condition {\em and} $\widetilde \nabla_a (G_{a,b}) = \widetilde \nabla_b (G_{a,b}).$
\end{tabular}
\end{defn}

\begin{figure} [h]
\centering {\includegraphics {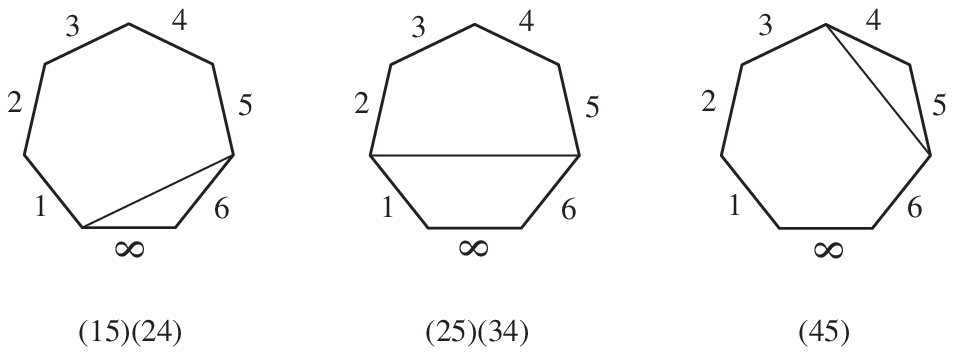}}
\caption{Examples of $\Gc \rightarrow \Sg_6$}
\label{f:mapphi}
\end{figure}

The machinery above is introduced in order to understand $\pi_1(\M{n})$. Fix an ordering of $\{1, 2, \ldots, n-1\}$ and use it to label the sides of each element in $\Gc$.  We define a map $\phi: \Cox{n-1} \rightarrow \Sg_{n-1}$ as follows:  Let $\phi(s_i)$ be the product of transpositions corresponding to the permuted labels of $G_i$ under $\widetilde \nabla_d (G_i)$. Figure~\ref{f:mapphi} gives a few examples.

It is not too difficult to show that the relations of $\Cox{n-1}$ carry over to $\Sg_{n-1}$.  Furthermore, the transpositions form a set of generators for $\Sg_{n-1}$, showing $\phi$ to be surjective.\footnote{To see this, simply consider the elements of $\Gc^2$.}  This leads to the following

\begin{thm} \textup{\cite[\S4]{djs2}}
\;$ker \, \phi \times \Z_2 \,=\, \pi_1 (\dM{n}) \times \Z_2 \,=\, \pi_1 (\M{n}).$
\end{thm}


\subsection{}
\begin{figure} [h]
\centering {\includegraphics {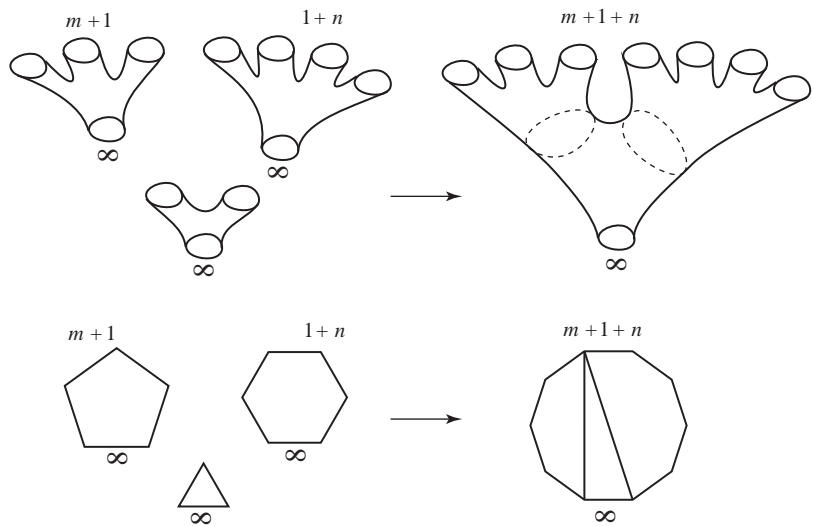}}
\caption{Pair-of-pants}
\label{pairpants}
\end{figure}

The {\em pair-of-pants} product (Figure~\ref{pairpants}) takes $m+1$ and $1+n$ marked points on $\R\Pj^1$ to $m+1+n$ marked points. The operad structure on the spaces \M{n+1}, its simplest case corresponding to the pair-of-pants product, defines composition maps \: $\Cox{m} \times \Cox{n} \rightarrow \Cox{m+n}$ \: analogous to the juxtaposition map of braids.

We can thus construct a monoidal category which has finite ordered sets as its objects and the group \Cox{n} as the automorphisms of a set of cardinality $n$, all other morphism sets being empty.  Note the following similarity between the braid group $\B_n$ obtained from the little cubes operad and the `quasibraids' \Cox{n} obtained from the mosaic operad:

\centerline{
\begin{tabular}{ccccc}
$\pi_1 (\C^n - \Delta)$ & $\rightarrowtail$ & $\B_n$
& $\twoheadrightarrow$ & $\Sg_n$ \\ [.4 cm]
$\pi_1 (\dM{n+1})$ & $\rightarrowtail$ & $\Cox{n}$ & $\twoheadrightarrow$ & $\Sg_n$
\end{tabular}} \medskip

\noindent There are deeper analogies between these structures which have yet to be studied.

%
%

\bibliographystyle{amsplain}

\end{document}